\documentclass[12pt]{article}
\usepackage{amsmath}
\usepackage{amssymb}
\usepackage{mathrsfs}
\usepackage{times}
\usepackage{epsfig,graphics}
\usepackage{color}

\usepackage{amsthm}

\begin{document}
\pagestyle{plain}
\hsize = 6.5 in
\vsize = 8.5 in
\hoffset = -0.5 in
\voffset = -0.5 in
\baselineskip = 0.29 in
\global\long\def\rd{{\rm d}}
\def\b {{\bf b}}
\def\vf {{\bf f}}
\def\vj {{\bf j}}
\def\vn {{\bf n}}
\def\vP{{\bf P}}
\def\vp{{\bf p}}
\def\x{{\bf x}}
\def\vu{{\bf u}}
\def\vv{{\bf v}}
\def\vw {{\bf w}}
\def\vz {{\bf z}}
\def\A{{\bf A}}
\def\B{{\bf B}}
\def\mC{{\bf C}}
\def\mD{{\bf D}}
\def\F{{\bf F}}
\def\mG{{\bf G}}
\def\I{{\bf I}}
\def\mM{{\bf M}}
\def\mS{{\bf S}}
\def\mQ{{\bf Q}}
\def\mR{{\bf R}}
\def\mT{{\bf T}}
\def\mU{{\bf U}}
\def\mV{{\bf V}}
\def\wtA{\widetilde{A}}
\def\mLa{\mbox{\boldmath$\Lambda$}}
\def\mGa{\mbox{\boldmath$\Gamma$}}
\def\mPhi{\mbox{\boldmath$\Phi$}}
\def\mPi{\mbox{\boldmath$\Pi$}}
\def\mXi{\mbox{\boldmath$\Xi$}}
\def\vg {\mbox{\boldmath$\gamma$}}
\def\vpi{\mbox{\boldmath$\pi$}}
\def\vphi{\mbox{\boldmath$\phi$}}
\def\vnu{\mbox{\boldmath$\nu$}}
\def\vkappa{\mbox{\boldmath$\kappa$}}
\def\vmu{\mbox{\boldmath$\mu$}}
\def\vxi{\mbox{\boldmath$\xi$}}
\def\vzeta{\mbox{\boldmath$\zeta$}}
\def\epr{\mbox{epr}}

\def\( {\left( }
\def\) {\right) }

\title{Stochastic Dynamics, Large Deviations Principle, and Nonequilibrium Thermodynamics}

\author{Liu Hong\\[5pt]
School of Mathematics,\\
Sun Yat-sen University, Guangzhou, 510275, P.R.C.\\
Email: hongliu@sysu.edu.cn\\[10pt]
and\\[10pt]
Hong Qian\\[5pt]
Department of Applied Mathematics,\\
University of Washington, Seattle, WA 98195-3925, U.S.A.\\
Email: hqian@uw.edu
}

\maketitle


\begin{abstract}
By examining the deterministic limit of a general $\epsilon$-dependent generator for Markovian dynamics, which includes the continuous Fokker-Planck equations and discrete chemical master equations as two special cases, the intrinsic connections among mesoscopic stochastic
dynamics, deterministic ODEs or PDEs, large deviations rate function, and macroscopic thermodynamic potential are established.
Our result not only solves the long-lasting question on the origin of entropy function in classical irreversible thermodynamics, but also reveals an emergent feature that arises automatically during the deterministic limit, through its large deviations rate function, with both time-reversible dynamics equipped with a Hamiltonian function and time-irreversible dynamics equipped with an entropy function.
\end{abstract}

\section{Introduction}

	Statistical equilibrium thermodynamics in terms of the theory of ensembles, as formulated originally by Gibbs, has a more fundamental origin; there is a growing consensus that the description of large deviations from the theory of probability provides a mathematical foundation for the subject \cite{ruelle-book}.  See \cite{barato-chetrite-2015,gq-pre-16,lu-qian-2020}, and the references within, for some of the recent developments.  One of the most important insights from \cite{gq-pre-16} is that one is able to derive nonequilibrium steady-state chemical thermodynamics using the same approach.

In a nutshell, the large deviations theory says that if a sequence of probability distributions has a deterministic limit, there exists concomitantly a scalar {\em rate function} for the convergence, which is defined on the entire state space.  This rate function provides a variational principle akin to those associated with various thermodynamic potentials for different ensembles \cite{Oono-89,Touchette-09,Smith-11,gq-ijmpb}.  The existence of a ``thermodynamic potential'' such as entropy, in addition to being justified through Clausius equality in traditional thermodynamics, and the logarithm of thermodynamic {\em Wahrscheinlichkeit} as Boltzmann did in the mechanical theory of heat, could be hidden in the {\em assumption} of the sequence of probability distributions and its deterministic limit.  We remind the readers that, according to Boltzmann's approach, the existence of a thermodynamic potential function, regardless of its specific form, is the only prerequisite for developing thermodynamics: Thermodynamic forces are introduced as conjugate variables, work as ``the force times the displacement'' then automatically becomes a part of the energy change.


It is well-known that the large deviations principle plays a key role in the equilibrium statistical mechanics \cite{ruelle-book,Ellis-06}.  In this study, we are going to show that it also acts as the statistical foundation of nonequilibrium thermodynamics, to be exact the Classical Irreversible Thermodynamics (CIT). CIT was developed by Onsager, Prigogine, de Groot and Mazur, \textit{et al.} since 1930s, as a natural generalization of concepts from equilibrium to nonequilibrium thermodynamics by introducing (i) the local equilibrium assumption which secures a spatiotemporal entropy density function and (ii) a linear thermodynamic force-flux relationships \cite{degroot-mazur}.  However, so far a solid and general theoretical foundation of CIT has not been established with a bottom-up approach. In this manuscript, we propose a bridge between the CIT and the application of large deviations principle (LDP) in the small-noise limit for a wide class of stochastic dynamics, representing by both diffusion and jump processes.

The existence of a nonequilibrium steady state is essential for our construction and its role cannot be overemphasized. In fact, a key step of our approach is to examine the stationary solution to a  Hamilton-Jacobi equation satisfied by the large deviations rate function (LDRF), which has been rigorously demonstrated for many cases in the limiting process of stochastic to deterministic nonlinear dynamics \cite{fengjin}. The latter is the macroscopic dynamics covered by CIT, for which the stationary large deviations rate function turns to be the (relative) entropy function.

To make our statement clear, the paper is organized as follows. In Section II, a very general $\epsilon$-dependent generator for Markovian dynamics, which includes the continuous Fokker-Planck equations and discrete chemical master equations as two special cases, is introduced. By examining the macroscopic limit, deterministic nonlinear ODEs and PDEs, which are widely used in modeling nonequilibrium processes, are derived from the Markovian dynamics as the noise level $\epsilon\rightarrow0$. Most importantly, during this limit process, stationary large deviations rate functions emerges automatically and serves as the entropy foundation for classical irreversible thermodynamics as shown in Section III. This conclusion is further justified through several concrete exactly-solvable examples.
Since the present paper is a synthesis of several lines of researches into a coherent theory for nonequilibrium thermodynamics, in Section IV, we discuss the relations between previous work and ours. The last section contains general discussions about entropy, energy, dissipative dynamics and conservative dynamics.

\section{Large deviations principle and stochastic dynamics}
\label{sec:2}

Most stochastic dynamics has a {\em natural} deterministic limit; which can and should be
understood as a form of the Law of Large Numbers (LLN) in the theory of probability.  If one uses an $\epsilon$ to represent this limiting process and denote the stochastic dynamics as
$Y_{\epsilon}(t)$, then one has
\begin{equation}
    \lim_{\epsilon\to 0} Y_{\epsilon}(t) = y(t),
\label{lln}
\end{equation}
where $y(t)$ represents a deterministic dynamics.  By envisioning the $y(t)$ as the limit of a sequence of $Y_{\epsilon}(t)$, there is an emergent variational principle associated with the deterministic $y(t)$: This is the essence of our new thermodynamics.  One example of the type of limit theorems in (\ref{lln}) is given by Kurtz. According to \cite{kurtz}, for any finite time the volume averaged trajectories of particle numbers inside a given chemical reaction system, which follows a pure jump Markov process, will converge in probability to the solution of a set of deterministic ordinary differential equations (ODEs) in the limit of infinitely large particle number and volume (with the volume per particle being $\epsilon$) while keeping a finite ratio between the two: the concentration.  One can find details of this result in the texts by van Kampen \cite{kampen}, Gardiner \cite{gardiner}, or Keizer \cite{keizer}.

In general $Y_{\epsilon}(t)$ can be either discrete or continuous: For Markov dynamics, it can be a diffusion process driven by a stochastic Wiener process $W(t)$, or an integer-valued, continuous-time jump process on the lattice $\mathbb{Z}^N$ driven by a stochastic Poisson process $\Xi(t)$.  In the former case, $\sqrt{\epsilon}$ is the amplitude of the $W(t)$; and for the latter, the jump size is
proportional to $\epsilon$, which means we taking the continuum limit. For a spatiotemporal dynamics
$Y_{\epsilon}(x,t)$, the well-known examples are the solution to a
stochastic partial differential equation (SPDE) with $\epsilon$ scaled white noise and interacting particle systems (IPS) \cite{liggett}.  The corresponding limit law as in (\ref{lln}) is called
the {\em hydrodynamic limit}: In this case, the limit $y(x,t)$ usually satisfies a nonlinear partial differential equation (PDE) \cite{gpv,htyau,fengjf}.

\subsection{$\epsilon$-dependent Markov transition probability}

Let us now consider the situation that a continuous-time
stochastic, Markov $Y_{\epsilon}(t)$ is not given explicitly, but
only in terms of a dynamic equation, as its generator:
\begin{subequations}
\begin{eqnarray}
    \frac{\partial}{\partial t} \rho_{\epsilon}(z;t) &=& \int_{\mathbb{R}^n}
             T(z|\xi;\epsilon) \rho_{\epsilon}(\xi;t) \rd\xi,
\label{ckeq}\\
	\rho_{\epsilon}(z;t)\rd z  &=&
             \Pr\big\{ z< Y_{\epsilon}(t) \le z+\rd z \big\}.
\end{eqnarray}
\end{subequations}
This is a special form of the Chapman-Kolmogorov equation
for Markov dynamics, assuming a time-homogeneous rate for the transition
probability. $\rho_{\epsilon}(z;t)$ denotes the probability density of the system in state $z$ at time $t$, with $\epsilon\ll 1$ as a small parameter indicating the strength of randomness.

$T(z|\xi;\epsilon)$ is the transition probability from state $\xi$ to state $z$, which has the following essential properties.  For
$z,\xi\in\mathbb{R}^n$ and $\epsilon > 0$,

(i)
\[
                T(z|\xi;\epsilon) \ge 0 \ \text{for} \ z\neq \xi;
\]

(ii)
\[           \int_{\mathbb{R}^n} T(z|\xi;\epsilon) \rd z = 0;
\]

(iii)
\[
          \lim_{\epsilon\to 0} T(z|\xi;\epsilon) = -F(\xi)\delta'(z-\xi).
\]
The first two are standard properties for an infinitesimal Markov generator, the third one indicates a deterministic (weak-noise) limit as $\epsilon\to 0$:
\begin{eqnarray}
     \frac{\partial \rho(z;t)}{\partial t}  &=& -\int_{\mathbb{R}^n}
           \Big(F(\xi)\delta'(z-\xi)\Big)  \rho(\xi;t) \rd\xi
\nonumber\\
	&=& -\int_{\mathbb{R}^n}
           \nabla\cdot\Big(F(\xi)\rho(\xi;t)\Big) \delta(z-\xi) \rd\xi
  \ = \ -\nabla\cdot\Big(F(z)\rho(z;t)\Big),
\end{eqnarray}
which represents a nonlinear ODE
\begin{equation}
     \frac{\rd z(t)}{\rd t} = F(z).
\label{the-de}
\end{equation}

Note that (i) and (ii) imply that $T(z|\xi;\epsilon)$ has a
negative signed Dirac-$\delta$ atomic measure at $z=\xi$.  However,
for finite $\epsilon$, there could be other positive signed atomic
measure at $z-\xi\neq 0$ in the form
\begin{equation}
  T(z|\xi;\epsilon) = R(\xi)
           \left[ \frac{ \delta(z-\xi+\epsilon\nu)-\delta(z-\xi) }{\epsilon}\right],
\end{equation}
in which the amplitude $R(\xi)$ is non-negative and location
$\nu$ can be positive or negative.  We therefore assume the following
general form for the Markov generator $T(z|\xi;\epsilon)$:\footnote{The theory of L\'{e}vy processes gives a more rigorous treatment of a Markov process whose probability of increments changes continuously in time \cite{cinlar-book,schilling-book}.  Our assumption here amounts to a class of diffusion with jumps in units of $\epsilon$.  }
\begin{subequations}
\label{T-expansion}
\begin{eqnarray}
  T(z|\xi;\epsilon) &=& \sum_{\ell=-m}^m  R_{\ell}(\xi)
        \left[ \frac{ \delta(\xi-z+\epsilon\nu_{\ell})-\delta(\xi-z)}{\epsilon} \right]
\\
     &-&  A(\xi)\delta'(z-\xi)  +
      \epsilon D(\xi)\delta''(z-\xi),
\end{eqnarray}
\end{subequations}
in which $A, \delta' \in\mathbb{R}^n$ are vectors and
$D, \delta'' \in\mathbb{R}^n\times\mathbb{R}^n$,
$D$ being a positive definite matrix; $\nu_{\ell}=-\nu_{-\ell}$. As a concrete example, it is noted that the term in (\ref{T-expansion}a) is the generator for chemical master equations with $\epsilon$ scaled jump size and time, leading to the $\epsilon\nu_{\ell}$ and $\epsilon^{-1}R_{\ell}$ \cite{dembo-book}, while the last two terms are those for the drift and $\epsilon$ scaled diffusion terms in the Fokker-Planck equation.

It is noted that in the above formula the mechanical and chemical contributions to the dynamics (respectively given by the $A$, $D$ and the $R$'s) scale with the same small parameter $\epsilon$. However, as the $\epsilon$ from chemical contribution originates from the infinitely large population limit; while the $\epsilon$ in the diffusion process represents thermal fluctuations in mechanical movements, they have not to be the same thing. This means in general we need to deal with a double-limit problem when studying its asymptotic behaviors (see footnote 1). Here, for simplicity, we just use the same small parameter $\epsilon$ for both processes, since in order to have a macroscopic chemomechanics they have to be infinitesimals on the same order.

In the limit of $\epsilon\to 0$, the generator in (\ref{T-expansion}) gives
\begin{equation}
  F(z) = A(z) + \sum_{\ell=-m}^m \nu_{\ell}R_{\ell}(z).
\label{eq-8}
\end{equation}
This provides a unified treatment of weak-noise limit of continuous diffusion processes as well as Kurtz's limit of jump process. The latter is modelled  by a Poisson process $Y(t;\lambda)$ with rate $\lambda$ to represent the particle number change during each chemical reaction in a stochastic way.
As $\epsilon\to 0$, $\epsilon Y(t; \epsilon^{-1}\lambda)\to \lambda t$.

\begin{figure}[t]
\[
\includegraphics[scale=0.55]{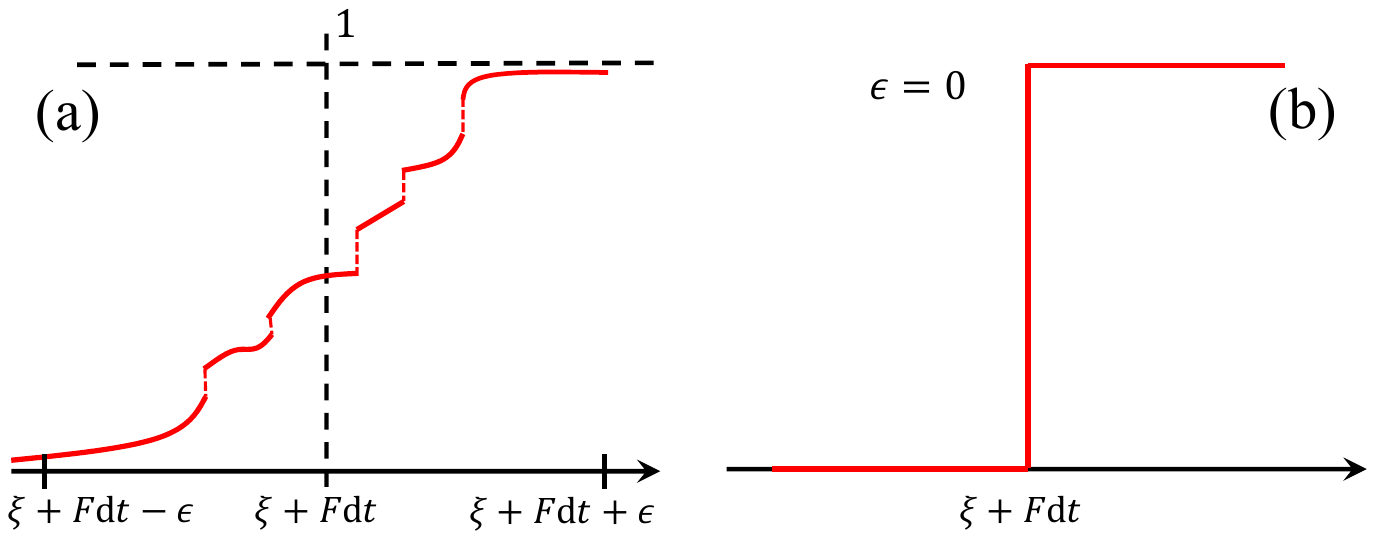}
\]
\caption{General, discontinuous transition probability distribution function
$P_{\epsilon}(x|\xi;\rd t)$ defined in (\ref{pepsilon}), shown in (a),
approaches to $H[x-\xi-F(\xi)\rd t]$, shown in (b), as $\epsilon\to 0$.
}
\label{fig-1}
\end{figure}

	Fig. \ref{fig-1} illustrates graphically the nature of the assumption
in (\ref{T-expansion}) when $z,\xi\in\mathbb{R}$:  The corresponding
transition probability distribution function in an
infinitesimal $\rd t$ is
\begin{equation}
          P_{\epsilon}(x|\xi;\rd t) = \int_{-\infty}^x
              \Big[ \delta(z-\xi)+T(z|\xi;\epsilon)\rd t \Big] \rd z.
\label{pepsilon}
\end{equation}
It in general contains discontinuous jumps.  In the limit of $\epsilon\to 0$, it converges to the Heaviside-step function $H\big[x-\xi-B(\xi)\rd t\big]$.

\subsection{Large deviations principle and Hamiltonian dynamics}

Besides the LLN, it is well-known that associated with the probability of $Y_{\epsilon}$ at time $t$, there exists a large deviations principle
\begin{equation}
\label{ldp}
    - \lim_{\epsilon\to 0} \Big(  \epsilon  \ln \Pr\big\{ \eta <
           Y_{\epsilon}(t) \le \eta+\rd\eta \big\} \Big)
    = \varphi(\eta,t),
\end{equation}
in which non-negative $\varphi(\eta,t)$, called the large deviations rate function, has a global
minimum zero when $\eta=y(t)$.

Now if one uses the result in (\ref{ldp}) as the basis for an
assumption like the WKB ansatz:
\begin{equation}
     \rho_{\epsilon}(z;t) = \exp\left(-\frac{\varphi(z,t)}{\epsilon}\right),
\end{equation}
and substitutes this expression into (\ref{ckeq}), one has the leading
order terms \footnotemark[1]:
\begin{eqnarray}
\frac{\partial\varphi(z,t)}{\partial t} &=& \epsilon \int_{\mathbb{R}^n}
      e^{\frac{\varphi(z,t)-\varphi(\xi,t)}{\epsilon}}  \Big\{
 A(\xi) \delta'(z-\xi)  -
      \epsilon D(\xi)\delta''(z-\xi)
\nonumber\\
	&-&  \left. \sum_{\ell=1}^m  R_{\ell}(\xi)
        \left[ \frac{ \delta(\xi-z+\epsilon\nu_{\ell})-\delta(\xi-z)}{\epsilon} \right]
           \right\} \rd\xi
\nonumber\\
	&=&  -A(z)\nabla\varphi(z,t)-\nabla\varphi(z,t)D(z)\nabla\varphi(z,t)
            -\sum_{\ell=-m}^m R_{\ell}(z) \Big[
              e^{\nu_{\ell}\nabla\varphi(z,t)} -1 \Big].  \hspace{1cm}
\label{hje}
\end{eqnarray}

\footnotetext[1]{Notice the subtlety for the following double limit as $\epsilon\to 0$
and $\epsilon'\to 0$, which is singular:
\begin{eqnarray*}
   && \lim_{\epsilon'\to 0} \lim_{\epsilon\to 0}\
          \epsilon' e^{\frac{\varphi(z)}{\epsilon'}} \int_{\mathbb{R}}
   e^{-\frac{\varphi(\xi)}{\epsilon'}}\left(
  \frac{  \delta(\xi-z-\epsilon\nu)-\delta(\xi-z)  }{\epsilon} \right) \rd\xi
\nonumber\\
  &=& \lim_{\epsilon'\to 0} \lim_{\epsilon\to 0}\  \frac{\epsilon'}{\epsilon}
        \Big[ e^{-\frac{\varphi(z+\epsilon\nu)-\varphi(z)}{\epsilon'}} -
    1  \Big]
\nonumber\\
	&=&  \left\{ \begin{array}{lcc}
                  \displaystyle   \lim_{\epsilon\to 0}\
        e^{-\frac{\varphi(z+\epsilon\nu)-\varphi(z)}{\epsilon}} - 1
        \ = \   e^{ -\nu[\rd\varphi(z)/\rd z] } - 1  &&  \epsilon'=\epsilon \\[12pt]
             \displaystyle   \lim_{\epsilon'\to 0}\
       - \frac{ \nu\epsilon'}{\epsilon'}
            \left[\frac{\rd\varphi(z)}{\rd z}\right]
        \ = \     -\nu \left[\frac{\rd\varphi(z)}{\rd z}\right]  &&
          \epsilon\to 0 \text{ first}  \end{array}
       \right.
\end{eqnarray*}
Actually, the limit does not exist if taking $\epsilon'\to 0$ first.}

With respect to the Hamilton-Jacobi equation in (\ref{hje}), it becomes possible to introduce a Hamiltonian function
\begin{align}
H(z,y)=&A(z)y+y^TD(z)y
            +\sum_{\ell=-m}^m R_{\ell}(z) \Big[
              e^{\nu_{\ell}y} -1 \Big],
\end{align}
and the corresponding Hamiltonian dynamics
\begin{align}
&\frac{\rd z}{dt}=\frac{\partial H(z,y)}{\partial y}=A(z)+2D(z)y+\sum_{\ell=-m}^m \nu_{\ell}R_{\ell}(z)e^{\nu_{\ell}y},
\label{Hamz}\\
&\frac{\rd y}{dt}=-\frac{\partial H(z,y)}{\partial z}=-A'(z)y-y^TD'(z)y-\sum_{\ell=-m}^m R'_{\ell}(z) \Big[e^{\nu_{\ell}y} -1 \Big].
\label{Hamy}
\end{align}
Clearly, the above Hamiltonian system accepts the zero-noise dynamics (\ref{the-de}) and (\ref{eq-8}) as a solution when the momentum equals to zero ($y=0$). So that it corresponds to the relaxation dynamics towards an attractor. All other solutions to the Hamiltonian system with $y\neq0$ corresponds to rare events which are impossible in the deterministic dynamics and are only populated in fluctuations.

Meanwhile, according to classical mechanics, we can also construct a variational principle by considering the Lagrangian function
\begin{align}
L(z,\dot{z})=&[y\dot{z}-H(z,y)]_{y=y(z,\dot{z})},
\end{align}
in which $y$ as a function of $z$ and $\dot{z}$ is obtained from solving the implicit equation (\ref{Hamz}). Clearly, $H(z,y)$ and $L(z,\dot{z})$ are Legendre transforms of each other. In terms of $L(z,\dot{z})$, the conjugate variable $y=\partial L(z,\dot{z})/\partial \dot{z}$. Consequently, the most probable path in consistent with above Hamiltonian dynamics with given $z(0)$ and $z(T)$ are given through the lease action principle
\begin{align}
\min_{z(s)} \int_0^TL(z(s),\dot{z}(s))ds.
\end{align}

\subsection{$-\varphi^{ss}$ as an entropy functional}

	We now show that the stationary solution to Eq. (\ref{hje}) is an
entropy functional for the nonlinear differential equation (\ref{the-de}):
\begin{eqnarray}
	\frac{\rd}{\rd t}\varphi^{ss}\big( z(t)\big) &=&
       F(z)\cdot \nabla_z\varphi^{ss}(z)
\nonumber\\
	&=& \left( A(z) + \sum_{\ell=-m}^m \nu_{\ell}R_{\ell}(z) \right) \cdot
             \nabla_z\varphi^{ss}(z).
\label{eq-13}
\end{eqnarray}
where $\varphi^{ss}(z)$ satisfies
\begin{equation}
  A(z)\nabla\varphi^{ss}(z) + \nabla\varphi^{ss}(z)D(z)\nabla\varphi^{ss}(z)
            + \sum_{\ell=-m}^m R_{\ell}(z) \Big[
              e^{\nu_{\ell}\nabla\varphi^{ss}(z)} -1 \Big] = 0.
\label{eq-14}
\end{equation}
From (\ref{eq-14}), and using inequality $e^a-1-a\ge 0$ for all $a\in\mathbb{R}$,
\begin{eqnarray}
	-F(z)\nabla\varphi^{ss}(z) &=&  \nabla\varphi^{ss}(z)D(z)\nabla\varphi^{ss}(z)
            + \sum_{\ell=-m}^m R_{\ell}(z) \Big[
              e^{\nu_{\ell}\nabla\varphi^{ss}(z)} - 1 -
     \nu_{\ell}\nabla\varphi^{ss}(z)  \Big]
\nonumber\\	
        &\ge& 0.
\end{eqnarray}
Therefore, the PDE (\ref{hje}) is the equation one seeks
to define an entropy, thus a nonequilibrium thermodynamics
of the nonlinear system (\ref{the-de}).


It is important to identify $-\rd\varphi^{ss}/\rd t$ not as entropy production rate,
rather as the instantaneous rate of entropy change.  Then one has
\cite{gq-pre-16,qian-arxiv}
\begin{eqnarray}
  \underbrace{ -\frac{\rd}{\rd t} \varphi^{ss}\big[z(t)\big] }_{\text{entropy change}}
    &=&   -\left[A(z) + \sum_{\ell=-m}^m
                    \nu_{\ell}R_{\ell}(z) \right] \nabla_z\varphi^{ss}(z)
\nonumber\\
   &=&  \underbrace{    A(z)D^{-1}(z)A(z) + \sum_{\ell=1}^m
                    \Big(R_{\ell}(z)-R_{-\ell}(z)\Big)\ln\left(\frac{R_{\ell}(z)}{R_{-\ell}(z)}\right)
             }_{\text{ entropy production}}
\nonumber\\
	&-&  \underbrace{ \Big(A+D\nabla\varphi^{ss}\Big)D^{-1}\Big(A+D\nabla\varphi^{ss}\Big)
           }_{\text{ mechanical drive}}
      -  \underbrace{  \sum_{\ell=1}^m \Big(R_{\ell}-R_{-\ell}\Big)\ln\left[\frac{R_{\ell}}{R_{-\ell}}
         e^{\nu_{\ell}\nabla\varphi^{ss}} \right]  }_{\text{chemical drive}}
\nonumber\\
	 &+&  \underbrace{ \Big( A\nabla\varphi^{ss} + \nabla\varphi^{ss}D\nabla\varphi^{ss}\Big)
           }_{\text{chemomechanical exchange}}
\label{dphidt}
\end{eqnarray}
According to (\ref{eq-14}), the last term representing chemomechanical exchange can also be expressed as
\[
      \Big[ A(z) + \nabla\varphi^{ss}(z) D(z)\Big]\nabla\varphi^{ss}(z)
          = - \sum_{\ell=-m}^m R_{\ell}(z) \Big[
              e^{\nu_{\ell}\nabla\varphi^{ss}(z)} -1 \Big].
\]
By ``chemomechanics'', we mean the continuous variables describing mechanical movements and the discrete jump processes representing chemical reactions as discrete events.  From the standpoint of all the atoms as point masses in the molecules, the distinction between mechanics and chemistry disappears; the latter is an emergent phenomenon of a very complex particle ``diffusion'' in a potential force field as first elucidated by H. A. Kramers \cite{kramers}.

	Let us now consider the specific situation in which both
the mechanical and chemical parts are in detailed
balance on their own \cite{gq-pre-16,qian-arxiv}:
\begin{eqnarray}
	A(z)	&=& -D(z)\nabla_z U(z),
\\[5pt]
	 \ln \left(\frac{R_{\ell}(z)}{R_{-\ell}(z)}\right) &=&
  - \nu_{\ell}\nabla_z G(z),
\end{eqnarray}
in which $U(z)$ is a ``mechanical'' potential function, and $G(z)$
is the Gibbs function for a chemical part.  Substituting these
two {\em potential conditions} into (\ref{the-de}) and (\ref{eq-8}),
we have
\begin{equation}
	\frac{\rd z(t)}{\rd t} = - \left[ D\nabla_z U +\sum_{\ell=1}^m
                    2\nu_{\ell}\hat{R}_{\ell}
       \sinh\left(\frac{1}{2}\nu_{\ell}\nabla_z G \right)
  \right].
\end{equation}
Eq. (\ref{hje}) becomes
\begin{equation}
\frac{\partial\varphi(z,t)}{\partial t}  = \big(\nabla U-\nabla\varphi \big) D\nabla\varphi
            -\sum_{\ell=0}^m 2\hat{R}_{\ell}  \left[ \cosh\left(
              \nu_{\ell}\nabla\varphi-\frac{1}{2}\nu_{\ell}\nabla G \right)
            - \cosh\left(\frac{1}{2}\nu_{\ell}\nabla G \right)
 \right],
\end{equation}
in which $\hat{R}_{\ell}(z)=[R_{\ell}(z)R_{-\ell}(z)]^{\frac{1}{2}}$.
In a chemomechanic equilibrium, both thermomechanics described by the continuous variables and thermochemistry represented by the jump processes have to be in their respective equilibrium; and furthermore, the chemomechanical energy transduction has to be precisely balanced by mechanochemical energy transduction in a reverse process.  $U$ and $G$ are actually different views of the same invariant probability measure with detailed balance. As a consequence, a global chemomechanical equilibrium is reached when $U(z)=\varphi^{ss}(z)=G(z)$.

\subsection{PDEs as deterministic limits}

Our previous derivations can be formally extended to nonlinear PDEs, which are most widely used models in nonequilibrium thermodynamics. PDEs originally rose from treating fluid
dynamics in terms of Newtonian mechanics; the thermodynamics of
continuum thus naturally follows. There is a long tradition in the
physics of nonequilibrium thermodynamics in terms of continuum theory
\cite{degroot-mazur} at one hand, and in formulating mathematical theory
of PDEs in terms of ideas from thermodynamics \cite{serrin-78,serrin-86}
on the other hand.  Ultimately, with a statistical foundation of thermodynamics
in mind, a system of PDEs can be understood as the hydrodynamic limit
of a SPDE or IPS.  For example, nonequilibrium thermodynamics emerges in
the asymptotic limit via the law of large deviations in the interaction particle system \cite{bertini-07};
generalized Gibbsian chemical thermodynamics emerges
in the asymptotic limit of Kurtz's theorem \cite{kurtz}; and so does the
chemomechanics we outlined in Sec. \ref{sec:2} above.

To begin with, let us consider a spatiotemporal stochastic process $Y_{\epsilon}(x,t)$,
where $x\in\Omega\subset\mathbb{R}$, whose probability distribution for
the entire function of $x$ at a give time $t$,
\begin{equation}
           \Pr\Big\{ z(x) < Y_{\epsilon}(x,t) \le z(x) + \rd z(x) \Big\}   =  \rho_{\epsilon}\big[z(x);t\big]\rd z(x),
\end{equation}
is given by a Chapman-Kolmogorov equation
\begin{equation}
\label{CKeq}
\frac{\partial}{\partial t}\rho_{\epsilon}[z(x);t]=\int T_{\epsilon}\big[z(x)\big|\xi(x)\big]\rho_{\epsilon}[\xi(x);t] \rd\xi(x),\\
\end{equation}
in which $\rho_{\epsilon}[z(x);t]=\rho_{\epsilon}[z(x)](t)$ is a functional of
$z(x)$, and $z(x)$, $\xi(x)$, and $\rd\xi(x)$ are all in an appropriate function
space.

Similarly, the transition probability $T_{\epsilon}\big[z(x)\big|\xi(x)\big]$ is
assumed to take the following general form
\begin{eqnarray}
   T_{\epsilon}\big[z(x)\big|\xi(x)\big]  &=&
     \sum_{\ell=-m}^m R_{\ell}[\xi(x)]\bigg[\frac{\delta[\xi(x)-z(x)+\epsilon\nu_{\ell}(x)]-\delta[\xi(x)-z(x)]}{\epsilon}\bigg]\nonumber\\
&-& A[\xi(x)]\delta'[\xi(x)-z(x)]+\epsilon D[\xi(x)]\delta''[\xi(x)-z(x)],
\end{eqnarray}
where $R_{\ell}[\xi(x)]$, $A[\xi(x)]$, and $D[\xi(x)]$ are all functionals of $\xi(x)$. $\delta'[\xi(x)-z(x)]$ and $\delta''[\xi(x)-z(x)]$ represent functional (or variational) derivatives of Dirac-$\delta$ functions defined as
\begin{align}
\delta'[z]\equiv \frac{\delta}{\delta z}\bigg(\delta[z]\bigg)=\lim_{h\rightarrow0}\frac{\delta[z+hdz]-\delta[z]}{h}.
\end{align}
With respect to this generator, in the limit of $\epsilon\rightarrow0$, we have
\begin{equation}
\lim_{\epsilon\rightarrow0}T_{\epsilon}[z(x)|\xi(x)]=-F[\xi(x)]\delta'[z(x)-\xi(x)],
\end{equation}
with functional
\begin{equation}
            F[z(x)]=A[z(x)]+\sum_{\ell=-m}^m\nu_{\ell}(x)R_{\ell}[z(x)];
\end{equation}
therefore,
\begin{align}
\frac{\partial}{\partial t}\rho_{0}[z(x);t]=-\frac{\delta}{\delta z}\bigg(F[z(x)]\rho_{0}[z(x);t]\bigg).
\end{align}
This equation for the functional $\rho_{0}[z(x);t]$ actually represents a nonlinear partial differential
equation
\begin{align}
\label{npde}
                        \frac{\partial z(x,t)}{\partial t}=F[z(x,t)],
\end{align}
in which $F$ maps a function space to which $z(x,t)$ belongs to $\mathbb{R}$.

\subsection{Conservation dissipation formalism}
We shall now assume a more concrete form for $F$ up to the second-order
spatial derivative of $z(x,t)$, $F[z(x,t)]=F(z,z_x,z_{xx},x)$.  This corresponds
to a rather broad class of nonlinear PDEs discussed in literature
\begin{equation}
\frac{\partial z(x,t)}{\partial t}=F(z,z_x,z_{xx},x),
\label{g-pde}
\end{equation}
in which $z_x$ denotes $\partial z(x,t)/\partial x$.
A PDE in which the $F=-\nabla\cdot j(z_x,z,x) + s(x)$
is called a transport equation; $F$ without the $z_{xx}$ term
is called {\em hyperbolic}, and with $z_{xx}$ term is called
{\em parabolic}. $F$ without the $z$ is called the Hamilton-Jacobi equation (HJE).  When $F$ does not contain $z_{xx}$, the nonlinear
first-order PDE can be solved by the {\em method of
characteristics} \cite{carrier-pearson,kevorkian-pde,evans}; for
an HJE this method gives rise to a Hamiltonian system.

In particular, a group of first-order PDEs in a form of
\begin{equation}
  \partial_tz=-\sum_{j=1}^n \partial_{x_i}J_i(z)+Q(z),
\end{equation}
where
\begin{equation*}
z=\left(\begin{array}{c}z_s\\z_d\end{array}\right),\quad J_i(z)=\left(\begin{array}{c}f_i(z)\\g_i(z)\end{array}\right),\quad Q(z)=\left(\begin{array}{c}0\\q(z)\end{array}\right),
\end{equation*}
are of great interest in both nonequilibrium thermodynamics and mathematical physics. $z=z(x,t)$ is a set of independent variables used for characterizing the system under study. $J_i(z)$ are fluxes along the $x_i$ direction, while $q(z)\neq0$ represents nonzero source or sink terms. We notice terms in $Q(z)$ corresponding to $z$ are all zeros, which means variables $z_s$ obey local conservation laws. The theoretical significance of local conservation laws is that they point out which kind of variables (an integration of $z_s$ in the whole space) do not change with time. In contrast, since the source terms for $z_d$ are nonzero, the spacial integration of $z_d$ is no more constant. This gives a natural classification of $z_s$ and $z_d$ variables.

With respect to above equations, a very general mathematical formulation -- Conservation Dissipation Formalism \cite{zhu-15}, for modeling nonequilibrium processes was constructed, two key assumptions of which read
\begin{itemize}
\item There is a strictly concave smooth function $\eta=\eta(z)$, called entropy, such that $\eta_{zz}\cdot J_{iz}(z)$ is symmetric for each $i$ and for all $z=(z_s, z_d)$ under consideration;
\item There is a positive definite matrix $M(z)$, called dissipation matrix, such that $q(z)=M(z)\cdot\eta_{z_d}(z).$
\end{itemize}
The first assumption is the famous entropy condition for hyperbolic conservation laws due to Godunov \cite{Godunov-61}, Friedrichs and Lax \cite{Lax-71} and \textit{et al.}, which ensures the system is globally symmetrizable hyperbolic. Then the Poincare lemma implies that there is a function $K_i=K_i(z)$ such that
\begin{equation}
\eta_z\cdot J_{iz}=K_{iz}.\nonumber
\end{equation}
The second condition is a nonlinearization of the celebrated Onsager's reciprocal relation \cite{Onsager-31-1, Onsager-31-2}, which ensures that the states far away from equilibrium tend to
equilibrium in the long time.

Now it is easy to see that, for $F[z(x,t)]$ in (\ref{g-pde}) which meet the two requirements of CDF, $-\varphi^{ss}$ turns to be the desired entropy function. And its corresponding time evolution is given through the following balance equation
\begin{equation}
\frac{\partial \varphi^{ss}(z)}{\partial t}=\sum_{j=1}^n \partial_{x_i}K_i(z)-\sigma(z),
\label{ebl}
\end{equation}
where $K_i(z)$ is the entropy flux, and $\sigma(z)=\frac{\partial \varphi^{ss}(z)}{\partial t}M(z)\frac{\partial \varphi^{ss}(z)}{\partial z}\geq0$ is the entropy production rate. This result establishes an interesting connection among stochastic thermodynamics, large deviations rate function and macroscopic nonequilibrium thermodynamics. As the formulation for CDF looks similar to what we did for CIT in the next section, no more details will be presented here. Interesting readers may work it out by themselves.

\section{Large deviations principle and classical irreversible thermodynamics}

\subsection{The logical structure of CIT}

In contrast to the deterministic limit of stochastic processes discussed in previous sections, which provides a direct linkage between mesoscopic and macroscopic dynamics, there are also other schools, like classical irreversible thermodynamics, trying to derive the governing equations for macroscopic deterministic dynamics directly from a thermodynamic point of view. The mathematics of macroscopic classical irreversible
thermodynamics, as presented in \cite{degroot-mazur},
has a very elegant and clear logical structure, which we summarize
here.

	(i) First, one considers the macroscopic system is locally fully specified by several quantities, say
$\vu_i,\vv_i,\vw_i$, $i\in\mathfrak{S}$, which are called ``state variables''. For example, in classical hydrodynamics, the fluid density $\rho$, velocity $v$ and total energy $e$ are most often used ones. This assumption is generally referred to as the ``local equilibrium hypothesis'' in literature, which allows the application of concepts and methodology in equilibrium thermodynamics directly to nonequilibrium systems. In CIT, another remarkable feature of state variables is that each of them satisfies a system of conservation law of its own, which means $\vu(t)=\{u_i(t),i\in\mathfrak{S}\}$ follows
\begin{equation}
    \frac{\rd u_i(t)}{\rd t} = \sum_{j\in\mathfrak{S}}
                   \Big( J^{(u)}_{ji}-J^{(u)}_{ij}\Big),
\end{equation}
where $J_{ij}^{(u)}\ge 0$ is a one-way flux.

	(ii) The local
equilibrium hypothesis also guarantees the existence of a
local strictly convex entropy function $s_i\equiv S(u_i,v_i,w_i)$.  Then by
differential calculus (or Gibbs relation in thermodynamics) one has
\begin{subequations}
\label{app-e1}
\begin{eqnarray}
     \frac{\rd s_i}{\rd t} &=&
    \left(\frac{\partial s_i}{\partial u_i}\right)
 \left[\frac{\rd u_i(t)}{\rd t}\right] +
    \left(\frac{\partial s_i}{\partial v_i}\right)
 \left[\frac{\rd v_i(t)}{\rd t}\right] +
   \left(\frac{\partial s_i}{\partial w_i}\right)
 \left[\frac{\rd w_i(t)}{\rd t}\right]
\nonumber\\
      &=&  \sum_{\xi=u,v,w}
      \left(\frac{\partial s_i}{\partial \xi_i}\right)
       \sum_{j\in\mathfrak{S}}
        \Big( J^{(\xi)}_{ji}-J^{(\xi)}_{ij}\Big)
\\
	&=& \underbrace{\sum_{\xi=u,v,w}\sum_{j\in\mathfrak{S}}
              \frac{1}{2}\left[\frac{\partial s_i}{\partial \xi_i}-
          \frac{\partial s_j}{\partial \xi_j}\right]
       \Big( J^{(\xi)}_{ji}-J^{(\xi)}_{ij}\Big)  }_{\text{local entropy production rate = force $\times$ flux}}  +
      \underbrace{ \sum_{\xi=u,v,w}\sum_{j\in\mathfrak{S}}
           J^{(S,\xi)}_{ji}
          }_{\text{entropy exchange flux}},
\end{eqnarray}
in which the net entropy flux due to transport of $\xi$:
\begin{equation}
    J^{(S,\xi)}_{ij} =
    \frac{1}{2}\left[\frac{\partial s_i}{\partial \xi_i}+
          \frac{\partial s_j}{\partial \xi_j}\right]
       \Big( J^{(\xi)}_{ij}-J_{ji}^{(\xi)}\Big) =
                -J^{(S,\xi)}_{ji} ,
\end{equation}
\end{subequations}

	Eq. (\ref{app-e1}) has established a {\em local}
entropy balance law in the form given by (\ref{ebl}).
If $\xi$ represents energy, volume, or the concentration
of a chemical species, then $(\partial s_i/\partial \xi_i)
\equiv (\partial S/\partial\xi)_i$, $i\in\mathfrak{S}$,
will be $1/T_i$, $p_i/T_i$, and $\mu_i/T_i$ respectively, with
$T_i$, $p_i$ and $\mu_i$ being local temperature,
pressure, and chemical potential.  Then the
corresponding {\em thermodynamic forces} between
states $i$ and $j$ are $(1/T_i-1/T_j)$, $(p_i/T_i-p_j/T_j)$,
and $\mu_i/T_i-\mu_j/T_j$.

	(iii) Eq. (\ref{app-e1}b), which splits (\ref{app-e1}a) into a
symmetric and an antisymmetric terms, is actually a discrete
version of the integration by parts in differential calculus,
which is employed in the third step of CIT based
on continuous variables.  We recognize $\sum_{j\in\mathfrak{S}}\big(J_{ij}-J_{ji}\big)$ as a discrete analogue of a divergence term, in which $J_{ij}\ge 0$ is a one-way flux.
Thus we have
\begin{equation}
	  \left(\frac{\partial s_i}{\partial \xi_i}\right) \sum_{j\in\mathfrak{S}} \big(J_{ij}-J_{ji}\big)
   = -\sum_{j\in\mathfrak{S}} \bigg[\left(\frac{\partial s_j}{\partial \xi_j}\right)-\left(\frac{\partial s_i}{\partial \xi_i}\right)\bigg]J_{ji}
    + \sum_{j\in\mathfrak{S}} \Bigg[\left(\frac{\partial s_j}{\partial \xi_j}\right)J_{ji}-\left(\frac{\partial s_i}{\partial \xi_i}\right)J_{ij}\Bigg].
\label{int-by-p}
\end{equation}
The last term is again a divergence term, which is determined
by only boundary values when summed over a set of $i$'s.
Interestingly, if we identify $(s_j/\xi_j-s_i/\xi_i)$ as a
``thermodynamic force'' between states $i$ and $j$,
the term $(s_j/\xi_j-s_i/\xi_i)J_{ji}$ is not consistent with
Onsager's entropy production rate: According to his
theory a ``thermodynamics flux'' is the net flux
$(J_{ij}-J_{ji})$, not one-way flux \cite{hill-book}.
This observation suggests that instead of Eq. (\ref{int-by-p}),
we should make $J_{ij}$ and $J_{ji}$ symmetric:
\begin{eqnarray}
	  &&\left(\frac{\partial s_i}{\partial \xi_i}\right)\sum_{j\in\mathfrak{S}} \big(J_{ij}-J_{ji}\big)
 = \left(\frac{\partial s_i}{\partial \xi_i}\right)\sum_{j\in\mathfrak{S}}  \left[ \left(\frac{J_{ij}-J_{ji}}{2}\right) - \left(\frac{J_{ji}-J_{ij}}{2}\right) \right]
\nonumber\\
   &=& -\frac{1}{2}\sum_{j\in\mathfrak{S}} \bigg[\left(\frac{\partial s_j}{\partial \xi_j}\right)-\left(\frac{\partial s_i}{\partial \xi_i}\right)\bigg]
        \big(J_{ji}-J_{ij}\big)
    + \frac{1}{2}\sum_{j\in\mathfrak{S}}
     \bigg[\left(\frac{\partial s_i}{\partial \xi_i}\right)+\left(\frac{\partial s_j}{\partial \xi_j}\right)\bigg]\big(J_{ji}-J_{ij}\big).
\end{eqnarray}
This is precisely the Eq. (\ref{app-e1}b).

	(iv) Now the fourth step in CIT is to introduce a
thermodynamic force-flux relationship:
\begin{equation}
	 \left[\begin{array}{c}
      J_{ij}^{(u)}-J^{(u)}_{ji} \\
      J_{ij}^{(v)}-J^{(v)}_{ji} \\
       J_{ij}^{(w)}-J^{(w)}_{ji} \end{array}\right]
         = -  \mM_{ij}(\vu,\vv,\vw)\left[\begin{array}{c}
     \frac{\partial s_i}{\partial u_i}-
          \frac{\partial s_j}{\partial u_j} \\[4pt]
      \frac{\partial s_i}{\partial v_i}-
          \frac{\partial s_j}{\partial v_j} \\[4pt]
       \frac{\partial s_i}{\partial w_i}-
          \frac{\partial s_j}{\partial w_j} \end{array}\right],
\label{app-or}
\end{equation}
where $\mM_{ij}$ is a $3\times 3$ positive definite symmetric
matrix. With this assumption, the local entropy production rate in
(\ref{app-e1}) is strictly positive except all forces and fluxes
are zero. When $\mM_{ij}(\vu,\vv,\vw)$ is evaluated at
an equilibrium ($J_{ij}=J_{ji}$), it becomes a constant matrix. In that case, Eq. (\ref{app-or}) is called Onsager's near equilibrium
linear force-flux relationship, which can be derived from the
principle of detailed balance.

	The force-flux relation needs not to be linear. Another
well-known example is
\begin{equation}
 \left[\begin{array}{c}
      \ln(J_{ij}^{(u)}/J^{(u)}_{ji}) \\
       \ln(J_{ij}^{(v)}/J^{(v)}_{ji}) \\
       \ln(J_{ij}^{(w)}/J^{(w)}_{ji}) \end{array}\right]
         = -  \mM_{ij}\left[\begin{array}{c}
     \frac{\partial s_i}{\partial u_i}-
          \frac{\partial s_j}{\partial u_j} \\[4pt]
      \frac{\partial s_i}{\partial v_i}-
          \frac{\partial s_j}{\partial v_j} \\[4pt]
       \frac{\partial s_i}{\partial w_i}-
          \frac{\partial s_j}{\partial w_j} \end{array}\right].
\label{aff-flux}
\end{equation}
This is Gibbs' chemical affinity-flux relationship.  It implies
for each and every set of cyclic indices  $i_0,i_1,\cdots,i_n,
i_{n+1}=i_0$ in the state space $\mathfrak{S}$:
\begin{equation}
     \sum_{k=0}^n \mM^{-1}_{i_kk_{k+1}}
          \left[\begin{array}{c}
      \ln(J_{i_ki_{k+1}}^{(u)}/J^{(u)}_{i_{k+1}i_k}) \\
       \ln(J^{(v)}_{i_ki_{k+1}}/J^{(v)}_{i_{k+1}i_k}) \\
       \ln(J^{(w)}_{i_ki_{k+1}}/J^{(w)}_{i_{k+1}i_k}) \end{array}\right]
      = \left[ \begin{array}{c}  0 \\ 0 \\ 0 \end{array} \right],
\label{app-pindep}
\end{equation}
which is known as chemical detailed balance.  The three
zeros in (\ref{app-pindep}) implies there exist three potential
functions on the state space $\mathfrak{S}$.
With (\ref{aff-flux}) the local entropy production rate in
(\ref{app-e1}) is non-negative, and it is equal to zero
if and only if $J^{(\xi)}_{ij}=J^{(\xi)}_{ji}$ for all
$i,j\in\mathfrak{S}$ and $\xi=u,v,w$.

\subsection{CIT for master equations}

	Let us now follow the same steps (i) to (iv) for
a master equation which conserves the probability
\begin{equation}
  \frac{\rd p_i(t)}{\rd t} = \sum_{j\in\mathfrak{S}}
             \big( J_{ji}-J_{ij}\big), \
                J_{ij} = p_i(t)q_{ij} \ge0.
\label{master-eq}
\end{equation}
Introducing a local entropy function
$s_i= -p_i \ln p_i$.  Then
\begin{eqnarray}
  \frac{\rd s_i}{\rd t} &=& -\big(\ln p_i + 1\big)\frac{\rd p_i(t)}{\rd t}
\nonumber\\
	&=&  \frac{1}{2} \sum_{j\in\mathfrak{S}} \big(p_iq_{ij}-p_jq_{ji}\big)\ln\left(\frac{p_i}{p_j}\right) +
 \frac{1}{2} \sum_{j\in\mathfrak{S}} \big(p_iq_{ij}-p_jq_{ji}\big)
   \big[ \ln\big(p_ip_j\big) + 2\big].
\label{me-dsdt}
\end{eqnarray}
Now introducing an affinity-flux relationship
$\ln(p_i/p_j)=M_{ij}(p_iq_{ij}-p_jq_{ji})$, where
\begin{equation}
           M_{ij}(p_i,p_j) =\frac{\ln p_i-\ln p_j }{
           p_iq_{ij}-p_jq_{ji}}.
\end{equation}
It is easy to show that each element of $M$ is strictly positive if
and only if $q_{ij}=q_{ji}$.

	When $q_{ij}\neq q_{ji}$, the above simple entropy function
that is independent of $\{q_{ij}\}$ can no longer be a valid
choice.  Rather, a proper entropy function has to be informed
by the dynamics in (\ref{master-eq}).  One of the best known
examples is to consider the stationary probability distribution to
(\ref{master-eq}) $\{\pi_i\}$: $\tilde{s}_i=-p_i\ln(p_i/\pi_i)$.
This is the fundamental idea of free energy.   Instead of
(\ref{me-dsdt}) one then has
\begin{equation}
 \frac{\rd\tilde{s}_i(t)}{\rd t} = \frac{1}{2} \sum_{j\in\mathfrak{S}} \big(p_iq_{ij}-p_jq_{ji}\big)\ln\left(\frac{p_i\pi_j}{p_j\pi_i}\right) +
 \frac{1}{2} \sum_{j\in\mathfrak{S}} \big(p_iq_{ij}-p_jq_{ji}\big)
   \left[ \ln\left(\frac{p_ip_j}{\pi_i\pi_j}\right) + 2\right].
\end{equation}
One therefore has an affinity-flux relationship
$\ln(p_i\pi_j/p_j\pi_i)=\tilde{M}_{ij}(p_iq_{ij}-p_jq_{ji})$, with
\begin{equation}
  \tilde{M}_{ij}(p_i,p_j) = \frac{\ln(p_i/\pi_i)-\ln(p_j/\pi_j) }{
           p_iq_{ij}-p_jq_{ji}}.
\end{equation}
The matrix $\tilde{M}$ is symmetric and semi-positive definite if and only if the detailed balance condition holds
$\pi_iq_{ij}=\pi_j q_{ji}$ (see e.g. Ref. \cite{Peng-18} for a rigorous proof).

\subsection{$\varphi^{ss}$ as the statistical foundation of CIT}
The entropy function plays a key role during the formulation of CIT, however its origin is a mystery in macroscopic thermodynamics and CIT does not provide an answer to it. Interestingly, the large deviation rate function obtained from the limit process of mesoscopic stochastic dynamics turns out to be the desired entropy function for the macroscopic thermodynamic modeling, and thus it provides a solid statistical foundation for CIT.

To make this point clear, we start with the stationary large deviation function (or the free energy function in this case) and examine its full time derivative in accordance with CIT
\begin{align}
&\frac{\rd\varphi^{ss}[z(x,t)]}{\rd t}=\frac{\rd z(x,t)}{\rd t}\frac{\delta\varphi^{ss}[z(x,t)]}{\delta z}\nonumber\\
=&-\sum_{l=-m}^m R_l[z]\bigg[e^{\nu_l\delta\phi^{ss}[z]/\delta z}-1\bigg]-A[z]\frac{\delta\varphi^{ss}[z]}{\delta z}-\frac{\delta\varphi^{ss}[z]}{\delta z}D[z]\frac{\delta\varphi^{ss}[z]}{\delta z}+\frac{dz}{dt}\frac{\delta\varphi^{ss}[z]}{\delta z}\nonumber\\
=&-\sum_{l=-m}^m R_l[z]\bigg[e^{\nu_l\delta\varphi^{ss}[z]/\delta z}-\nu_l\frac{\delta\varphi^{ss}[z]}{\delta z}-1\bigg]+\bigg[\frac{dz}{dt}-\sum_{l=-m}^m \nu_l R_l[z]-A[z]-\frac{\delta\varphi^{ss}[z]}{\delta z}D[z]\bigg]\frac{\delta\varphi^{ss}[z]}{\delta z}\nonumber\\
=&-\sigma_1-\sigma_2.
\end{align}
It is seen that $\sigma_1\geq0$ by Bernoulli's inequality. While to keep $\sigma_2\geq2$ in accordance with the second law of thermodynamics, CIT suggests to take
\begin{align}
&\frac{\rd z}{\rd t}-\sum_{l=-m}^m \nu_l R_l[z]-A[z]=\frac{\delta\varphi^{ss}[z(x,t)]}{\delta z}\bigg(D[z]-M[z]\bigg),
\end{align}
where $M[z]\geq0$ must be semi-positive definite. In particular, if we choose $M[z]=D[z]$, the macroscopic equation in (\ref{npde}) is recovered. Comparing to the original equation, we see that models derived from CIT are not completely specified unless the entropy production rate is given too (which means $M[z]$ is given).  This ambiguity is arised from the fact that a dissipative process is not fully specified by the entropy function, but also by its dissipation rate.

\begin{figure}[t]
\[
\includegraphics[scale=0.55]{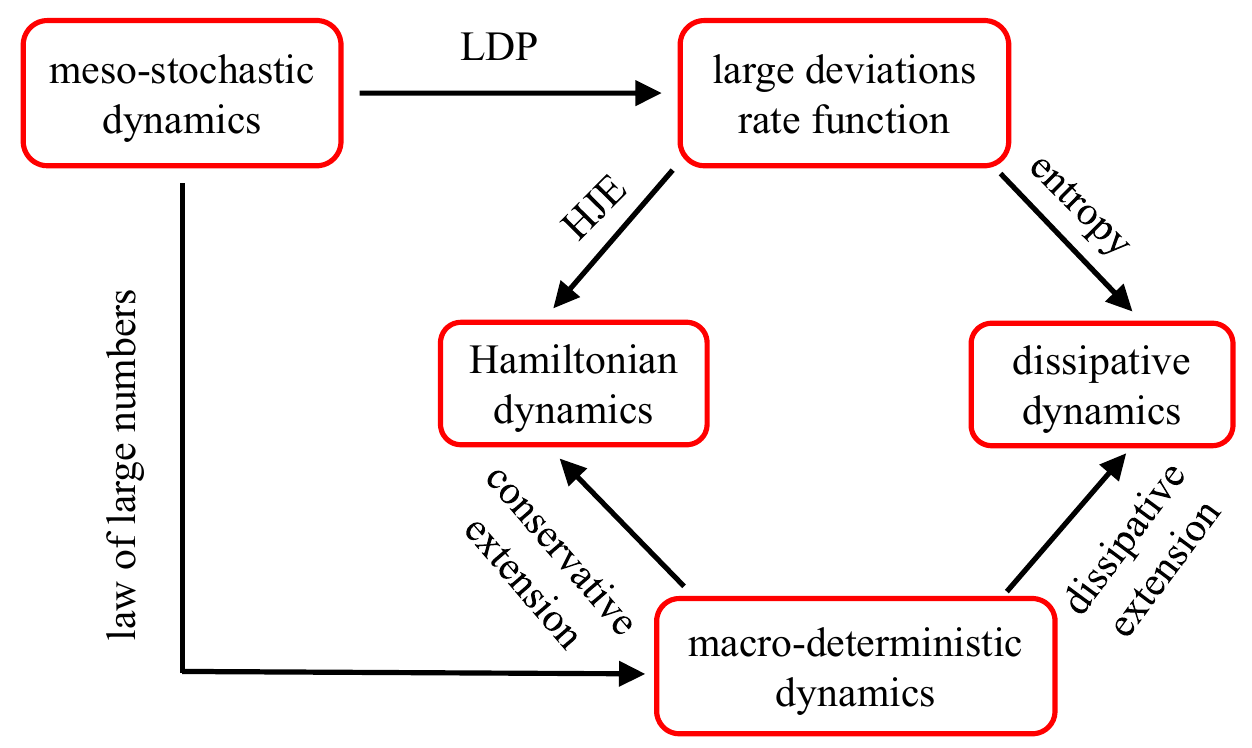}
\]
\caption{Relations among stochastic dynamics, macroscopic limit, large deviations theory, \textit{etc.} LDP: large deviations principle, HJE:
Hamilton-Jacobi equation.}
\label{fig-2}
\end{figure}

\subsection{Exactly-solvable models}

Finally, we look at several examples, which could be explicitly solved, to illustrate the intrinsic relations among mesoscopic stochastic dynamics, macroscopic deterministic dynamics, large deviations rate function, classical irreversible thermodynamics, Hamiltonian dynamics, and so on; see Fig. 2. \\

\subsubsection{Fokker-Planck equations for the Ornstein-Uhlenbeck process}

By taking $R(\xi)=0, A(\xi)=-az, D(\xi)=D$ in the generator in (\ref{T-expansion}), we arrive at the famous Fokker-Planck equation
\begin{eqnarray}
\frac{\partial p_{\epsilon}(z,t)}{\partial t}=
\frac{\partial}{\partial z}\cdot\bigg[\epsilon D\frac{\partial}{\partial z}p_{\epsilon}(z,t)+azp_{\epsilon}(z,t)\bigg].
\end{eqnarray}	
It corresponds to the Ornstein-Uhlenbeck process (OUP), a particular realization of the general Langevin dynamics, which reads
\begin{align}
\rd z(t)=-az\rd t+\sqrt{2\epsilon D}\rd B(t),
\end{align}
under the meaning of It\^{o}'s calculus. In this case, the distribution function could be exactly solved as
\begin{align}
p_{\epsilon}(z,t)=\bigg[\frac{a}{2\pi\epsilon D(1-e^{-2at})}\bigg]^{1/2}\exp\bigg[-\frac{az^2}{2\epsilon D(1-e^{-2at})}\bigg],
\end{align}
with respect to the initial condition $p_{\epsilon}(z,0)=\delta(z)$.

It is straightforward to show the large deviations rate function $\varphi(z,t)=az^2/[2D(1-e^{-2at})]$ and its stationary solution $\varphi^{ss}(z)=az^2/(2D)$. With respect to these formulas, we can repeat previous derivations of CIT. And it is easy to check that the relation $\rd z/\rd t=-az$ guarantees a positive entropy production. On the other hand, as suggested by the large deviations principle, we can also introduce a Hamiltonian dynamics
\begin{align}
&\frac{\rd z}{\rd t}=2Dy-az,\\
&\frac{\rd y}{\rd t}=ay,
\end{align}
with a Hamiltonian function $H(z,y)=Dy^2-azy$, which is equivalent to a Lagrangian dynamics
\begin{align}
\ddot{z}-a^2z=0
\end{align}
given by the Lagrangian function $L(z,\dot{z})=(\dot{z}+az)^2/(4D)$. Noticeably, both dynamics are time-reversible generalizations of
$\rd z/\rd t=-az$.

\subsubsection{Chemical reactions under complex balance condition}

In the next example, we consider a discrete generator with $A(\xi)=0, D(\xi)=0$. In this case, the chemical master equations are obtained, whose deterministic limit gives usual ODEs
\begin{align}
\frac{\rd z(t)}{\rd t}=\sum_{\ell=1}^m \nu_{\ell}\cdot\big[R_{\ell}(z)-R_{-\ell}(z)\big]
\label{chem}
\end{align}
for $m$ chemical reactions with general rate functions $R_{\ell}(z)$. $z=(z_1,z_2,\cdots,z_n)$ are the concentrations of the $n$ species, and stoichiometric coefficients $\nu_{\ell}=(\nu_{\ell1},\nu_{\ell2},\cdots,\nu_{\ell n})$.

Thanks to the condition of complex balance, a concept first introduced by Horn and Jackson in 1972 \cite{Horn-72}, for a class of chemical reactions with laws of mass action, it can be shown that the kinetics equation above has a unique stationary solution $z^{ss}$ \cite{Anderson-15}, and the stationary large deviations rate function \cite{gq-pre-16}
\begin{align}
\varphi^{ss}(z)=\sum_{i=1}^n z_i\ln\bigg(\frac{z_i}{z^{ss}_i}\bigg)-z_i+z^{ss}_i
\end{align}
is a solution to
\begin{align}
\sum_{\ell=1}^m R_{\ell}(z) \Big[e^{\nu_{\ell}\nabla\varphi^{ss}(z)} -1 \Big]+R_{-\ell}(z) \Big[e^{-\nu_{\ell}\nabla\varphi^{ss}(z)} -1 \Big]=0.
\end{align}
Then according to CIT, a possible dissipative extension of (\ref{chem}) is
\begin{align}
&\frac{\rd z}{\rd t}=\sum_{\ell=1}^m \nu_{\ell}\cdot\big[R_{\ell}(z)-R_{-\ell}(z)\big]-M(z)\ln\bigg(\frac{z_i}{z^{ss}_i}\bigg),
\end{align}
which, however, requires a preknowledge of the stationary solution $z^{ss}$ and is nearly impossible for real applications. In contrast, the conservative extension does not require such kind of information. With respect to the Hamiltonian function $H(z,y)=\sum_{\ell=1}^m \big\{R_{\ell}(z) \big[
              e^{\nu_{\ell}y} -1 \big]+R_{-\ell}(z) \big[
              e^{-\nu_{\ell}y} -1 \big]\big\}$,
the derivation of corresponding Hamiltonian dynamics is straightforward, \textit{i.e.}
\begin{align}
&\frac{\rd z}{\rd t}=\sum_{\ell=1}^m \nu_{\ell}\big[R_{\ell}(z)e^{\nu_{\ell}y}-R_{-\ell}(z)e^{-\nu_{\ell}y}\big],\\
&\frac{\rd y}{\rd t}=-\sum_{\ell=1}^m \bigg\{R'_{\ell}(z) \Big[e^{\nu_{\ell}y} -1 \Big]+R'_{-\ell}(z) \Big[e^{-\nu_{\ell}y} -1 \Big]\bigg\}.
\end{align}
It is noted that when the momentum $y=0$, we recover original kinetic equations in (\ref{chem}), which is in fact dissipative and time irreversible in nature.

\subsubsection{Slow chemomechanical coupling near equilibrium}

Compared to purely diffusive or purely chemical processes, the non-trivial chemomechanical coupling is far more interesting.  In this case, we need to solve the stationary Hamilton-Jacobi equation in the full form,
\begin{equation}
      \Big[ A(z) + \nabla\varphi^{ss}(z) D(z)\Big]\nabla\varphi^{ss}(z)
          = \sum_{\ell=1}^m R_{\ell}(z) \Big[1-
              e^{\nu_{\ell}\nabla\varphi^{ss}(z)} \Big]
  +R_{-\ell}(z) \Big[1-
              e^{-\nu_{\ell}\nabla\varphi^{ss}(z)} \Big],
\label{eq62}
\end{equation}
which for most situations can only be solved numerically.

	We now consider the problem in which both the mechanical and chemical parts are in rapid equilibrium, e.g., detailed balance, on their own:
\begin{eqnarray*}
	A(z)	&=& -D(z)\nabla U(z),
\\[5pt]
	 \ln \left(\frac{R_{\ell}(z)}{R_{-\ell}(z)}\right) &=&
  - \nu_{\ell}\nabla G(z),
\end{eqnarray*}
but the chemomechanical conversion is slow and is not yet in equilibrium; thus $U(z)\neq G(z)$; the $\varphi^{ss}(z)$ then is equal to neither.  Substituting these two into (\ref{eq62}), we have
\begin{equation}
   \nabla\varphi^{ss}(z) D(z) \nabla\Big(\varphi^{ss}(z)-U(z)\Big)
  = \sum_{\ell=1}^m \Big[e^{\nu_{\ell}\nabla\varphi^{ss}(z)}
     -1 \Big] R_{\ell}(z) \left[
       e^{\nu_{\ell}\nabla \big( G(z)-\varphi^{ss}(z) \big) }-1 \right].
\label{eq63}
\end{equation}
The terms $\nabla(\varphi^{ss}-U)$ and $\nabla(\varphi^{ss}-G)$ are
thermomechanical force and thermochemical force respectively.  A true equilibrium has both being zero.  We see that if $\nabla(U-\varphi^{ss})>0$, then $\nabla(G-\varphi^{ss})<0$. This implies a net mechanical to chemical energy conversion.

Near equilibrium, the last factor on the right-hand side of (\ref{eq63}) can be approximated by Taylor expansion. In the special case of $D(z)=\sum_{\ell=1}^m \frac{R_{\ell}(z)-R_{-\ell}(z)}{\ln R_{\ell}(z)-\ln R_{-\ell}(z)}\nu_{\ell}\bigotimes\nu_{\ell}$ (called biochemical conductance in stoichiometric network theory), where $\bigotimes$ denotes the direct product of vectors, Eq. (\ref{eq63}) can be solved explicitly. In this case, the thermomechanical force and thermochemical force are equal and opposite, and are given by the difference between chemical potential and mechanical potential, i.e. $\nabla(\varphi^{ss}-U)=-\nabla(\varphi^{ss}-G)=\nabla(G-U)/2$.

For chemical reactions with the mass-action law in equilibrium, an explicit formula for $G(z)$ is known, that is $G(z)=\sum_{i=1}^n$ $z_i\ln(z_i/z^{ss}_i)-z_i+z^{ss}_i$ (see the second example above). As a consequence, $U(z)=\sum_{i=1}^n z_i\ln(z_i/z^{ss}_i)$ $-z_i+z^{ss}_i$ and $A(z)=-D(z)\cdot\ln(z/z^{ss})$. It is worth noting that the term $\sum_{i=1}^n z_i\ln z_i$, which appears in both $G(z)$ and $U(z)$, has different  interpretations in chemistry and in mechanics: The former is caused by chemical affinities, while the latter, according to Flory and Huggins, is an entropic effect reflecting the tendency of particle mixing.

\section{Relations to previous work}

There are three lines of researches that are highly relevant, and
with respect to which the current work is seeking a synthesis.

(i) The investigations of statistical behavior of a stochastic system, its fluctuations and its entropy productions.  The classical theory of Einstein, Landau, Onsager-Machlup, {\em etc.} for equilibrium fluctuations with detailed balance,
Keizer's for nonequilibrium steady-state fluctuations \cite{keizer}, and recently developed stochastic thermodynamics of finite trajectories \cite{Seifert-12} are well-known landmarks.  See
\cite{maes} for a recent development.

(ii) The study of large deviations beyond (i) for a {\em sequence} of stochastic systems with a LLN; the focuses here are to secure the existence of a rate function, and to derive its particular form.  Mathematical work tends to focus on the former while physics
literature on the latter.  For the mathematical method, we refer interested readers to the comprehensive reviews written for physicists \cite{Oono-89,Touchette-09,Smith-11,gq-ijmpb,Ellis-06}.
The traditional Gaussian fluctuation theory is simply the local Hessian structure near the global minimum of the rate function.

 	Along this direction, the weak-noise limit of Markovian dynamics, {\em e.g.} $\epsilon\to 0$ discussed in Sec. 2, has been discussed many times in the past, both rigorously and applied.  We shall particularly mention the extensive studies carried out by Graham and T\'{e}l \cite{Graham-84, Graham-85} for the Fokker-Planck Equation, and by Hu \cite{Hu-87}, Dykman {\em et. al.} \cite{Dykman-94}, and Agazzi {\em et. al.} \cite{agazzi} for the Chemical Master Equation.  

	(iii)  The parallel work for spatiotemporal, infinite-dimensional, stochastic system is technically much more demanding.  The demonstration of a deterministic limit and its form, usually a nonlinear PDE, is already a challenging task.   There is a large literature on obtaining the hydrodynamic limit. Its history, since 1980s, dates back to the pioneering works of Liggett \cite{liggett}, Guo {\em et. al.} \cite{gpv}, and many others.  On discrete IPS, the exact results on Asymmetric Simple Exclusion Process is particularly worth mentioning \cite{derrida-pp-98}. See \cite{Francoa-10} for a more recent study.
A very related approach is the reaction diffusion master equation (RDME) in continuous time:  See the earlier work of Keizer \cite{keizer} on the fluctuating Boltzmann equation, Vance and Ross \cite{vance-ross} on fluctuating Turing patterns, and the more recent \cite{hellander} in connection to numerical computations.  On continuous space-time SPDE, Graham and T\'{e}l also investigated Ginzburg-Landau equation with weak noise \cite{graham-tel-GLE}; Gon\c{c}alves and Jara have studied the nonequilibrium fluctuations of Kardar-Parisi-Zhang equation in terms of an Einstein relation \cite{goncalves}. With the LLN in hand, a systematic treatment of the fluctuations of hydrodynamic equations was carried out in Macroscopic Fluctuation Theory \cite{Bertini-02,Bertini-15}. Results on large deviations followed \cite{Kipnis-89,Derrida-98,bertini-07}.

	In contradistinction to the abovementioned (i)-(iii), our present work is a part of the recent development on setting the large-deviation structure, being a limiting law for a {\em sequence} of Markov processes, as the mathematical foundation for  nonequilibrium stochastic thermodynamics, beyond the studies of a Markov process.  This type of limit laws is to thermodynamic behavior what the central limit theorem is to the Gaussian fluctuation theory.  The focus is on the mathematical origin of thermodynamic behavior itself: To our best knowledge, it is the first time to combine all these results together, by using an $\epsilon$-dependent Markov generator, in providing a unified mathematical physics in which the LDT serves the statistical foundation of general nonequilibrium thermodynamics, like CIT and CDF.  This synthesis is missing from most of the previous studies. One exception is the Macroscopic Fluctuation Theory, developed by Jona-Lasinio and coworkers, in which thermodynamic relationships among force, work, and quasi-potential as energy were discussed; another is a study from us on extended irreversible thermodynamics \cite{Hong-20}.

	Because of the nature of synthesis, in Sec. 2.4 we discuss how to derive PDEs as the deterministic limits from a spatiotemporal stochastic process, in order to incorporate the general theory of CDF which includes both time and space.  The work in (iii) above provides the more advanced, in depth materials for this section. We merely give a heuristic coverage before moving toward the main purpose of our current work: Go beyond the law of large numbers, e.g., the proper hydrodynamic limit, and focus on the entropy structure on top of those hydrodynamic equations.

\section{Conclusion}
The macroscopic limit of mesoscopic stochastic dynamics, especially the Markovian dynamics with either continuous or discrete state space, is well understood since the pioneering works of Kurtz, Guo \textit{et. al.}, and many others.  On the other hand, the fact that the large deviations principle, which emerges concomitantly during the limit process and provides the entropy as a macroscopic potential function of nonequilibrium thermodynamics, was not fully appreciated in the past. In the current study, by examining the deterministic limit of a general $\epsilon$-dependent generator for Markovian dynamics, which includes the continuous Fokker-Planck equations and discrete chemical master equations as two special cases, the intrinsic connections among mesoscopic stochastic dynamics, its macroscopic limit, large deviations rate function, classical irreversible thermodynamics and its potential are established. To provide concrete examples on our construction, the purely linear Ornstein-Uhlenbeck process, chemical reactions under complex balance condition, and non-trivial slow chemomechanical coupling  near equilibrium are solved explicitly. The investigation of more interesting chemical-mechanical coupled systems is left to future studies. In conclusion, our result not only solves the long-lasting question on the origin of entropy function in CIT, and also suggests a more general principle for emergent phenomena.

In our study, an amazing observation is that both the dissipative dynamics equipped with a (relative) entropy function and the conservative dynamics equipped with a Hamiltonian function arise automatically from the large deviations principle of mesoscopic stochastic dynamics. This emergent phenomena not only highlights the inseparable nature of the first law and the second law of thermodynamics, which state the essential roles of energy and entropy in a thermodynamical view of dynamics, but also provides a practical way for constructing either conservative or dissipative dynamics of any given deterministic dynamics by considering its stochastic correspondence.

\section*{Acknowledgements}
L. H. acknowledges the financial supports from the National Natural Science Foundation of China (Grant \# 21877070) and the Hundred-Talent Program of Sun Yat-Sen University. H. Q. is partially
supported by the Olga Jung Wan Endowed Professorship.

\end{document}